 \documentclass[11pt]{article}
\newcommand{\version}{July 11, 2005}

\usepackage{amsthm,amsfonts,amsmath}
\swapnumbers
\theoremstyle{plain}

\newtheorem{thm}{THEOREM}[section]
\newtheorem{lm}[thm]{LEMMA}
\newtheorem{cl}[thm]{COROLLARY}

\theoremstyle{definition}

\theoremstyle{remark}
\newcommand{\upchi}{\raise1pt\hbox{$\chi$}}
\newcommand{\R}{{\mathord{\mathbb R}}}

\newcommand{\C}{{\mathord{\mathbb C}}}

\renewcommand{\|}{{\Vert}}
\numberwithin{equation}{section}
\begin{document}

%%%%%%%%%%%%%%%%%%%%%%
\markboth{\scriptsize{CLL \version}}{\scriptsize{CLL \version}}

\def\sn{{{\cal S}^N}}

\title{\bf{AN INEQUALITY OF HADAMARD TYPE FOR PERMANENTS}}
\author{\vspace{5pt} Eric Carlen$^1$, Michael Loss$^1$
Elliott H. Lieb$^{2}$ \\
\vspace{5pt}\small{$1.$ School of Mathematics, Georgia Tech,
Atlanta, GA 30332}  \\
\vspace{5pt}\small{$2.$ Departments of Mathematics and Physics, Jadwin
Hall,} \\[-6pt]
\small{Princeton University, P.~O.~Box 708, Princeton, NJ
  08544}\\
 }
\date{\version}
\maketitle
\footnotetext                                                                         
[1]{Work partially
supported by U.S. National Science Foundation
grant DMS 03-00349.    }                                                          
\footnotetext
[2]{Work partially
supported by U.S. National Science Foundation
grant PHY 01-39984.\\
\copyright\, 2005 by the authors. This paper may be reproduced, in its
entirety, for non-commercial purposes.}
                                                                      
\begin{abstract}

Let $F$ be an $N\times N$ complex matrix whose $j$th column is
the vector $\vec f_j$ in $\C^N$.  Let $|\vec f_j|^2$ denote the
sum of the absolute squares of the entries of $\vec f_j$. Hadamard's
inequality for determinants states that $|\det(F)| \le \prod_{j=1}^N|\vec
f_j|$. Here we prove a sharp upper bound on the permanent of $F$, which is
${\displaystyle |{\rm perm}(F)| \le {N!\over N^{N/2}} \prod_{j=1}^N|\vec
f_j|}$, and we determine all of the cases of equality.

We also discuss the case in which $|\vec f_j|$ is replaced by the
$\ell_p$ norm of the vector $\vec f$ considered as a function on
$\{1,2,\dots,N\}$. We note a simple sharp inequality for $p=1$, and
obtain bounds for intermediate $p$ by interpolation. The interpolated
bounds are not sharp, though there is a natural conjecture for what the
sharp bounds should be.

\end{abstract}

\medskip
\centerline{Mathemematics subject classification numbers:  15A45, 49M20}

\section{Introduction} \label{intro}

Let $F$ be an $N\times N$ complex matrix whose $j$th column is the
vector $\vec f_j$ in $\C^N$.  Let $|f_j|^2$ denote the sum of the
absolute squares of the entries of $\vec f_j$. Hadamard's inequality
for determinants \cite{H} states that $|\det(F)| \le \prod_{j=1}^N|\vec
f_j|$. Here we prove a sharp upper bound on the permanent of $F$:

\begin{thm} \label{thm1}
 For any  vectors $\vec f_1, \dots, \vec f_N$ in $\C^n$ we have the inequality
\begin{equation}\label{main}
|{\rm perm}(F)| \le {N!\over N^{N/2}}\prod_{j=1}^N|\vec f_j|\ .
\end{equation}
For $M>2$,  there is equality in (\ref{main})
if and only if at least one of the vectors $\vec f_j$ is zero, or else  $F$ is a rank one matrix and moreover,
each of the vectors $\vec f_j$ is a constant modulus vector; i.e., its entries all have the same absolute value.
\end{thm}

The conditions for equality can be reformulated as follows: There is equality in (\ref{main}) if and only if
one or more of the vectors $\vec f_j$ is zero, or else there are numbers $r_j$, $\xi_j$ $\zeta_j$, $j= 1,\dots,N$,
with each $r_j>0$ and each $\xi_j$ and $\zeta_j$ lying on the unit circle in the complex plane,
so that 
$$F_{j,k} =\xi_j\zeta_k r_k$$
for each $j,k$.

We shall give two proofs of this inequality. The first turns on recognizing (\ref{main}) as a close relative of the Brascamp--Lieb type inequality that we recently proved \cite{CLL} for integrals of products of functions on the sphere $S^N$.  To explain this way of viewing (\ref{main}), we first introduce some notation and terminology.

Let $\sn$ denote the symmetric group on $N$ letters; i.e., the group of all permutations $\sigma$ of
$\{1,\dots,N\}$.  Let the (composition) product in $\sn$ be denoted by juxtaposition, and for each
$1\le i,j\le N$ with $i\ne j$, let $\sigma_{i,j}$ be the pair permutation with $\sigma_{i,j}(i)=j$, $\sigma_{i,j}(j)=i$,
$\sigma_{i,j}(k) = k$ for $k\ne i,j$.
Let $\mu$ denote the uniform probability measure on $\sn$ so that if $g$ is any function on $\sn$,
\begin{equation}\label{unifrom}
\int_\sn g(\sigma){\rm d}\mu = {1\over N!}\sum_{\sigma\in\sn}g(\sigma)\ .
\end{equation}

 We may identify vectors in $\C^N$ with complex valued functions on  $\{1,\dots,N\}$ as follows:
If $f:\{1,\dots,N\}\to \C$, let $\vec f$ be the vector in $\C^N$ whose $j$th entry is $f(j)$. Conversely, given
a vector $\vec f$ in $\C^N$, define the function $f$ by setting $f(j)$ equal to the $j$th component of  $\vec f$.

For $1\le k\le N$, define the function $\pi_k:\sn\to \{1,\dots,N\}$ by
$$\pi_k(\sigma) = \sigma(k)\ .$$  If  $f:\{1,\dots,N\}\to \C$, then
$$f\circ \pi_k:\sn\to \C\ .$$

Let $\{f_1,\dots,f_N\}$ be any $N$ complex valued functions on  $\{1,\dots,N\}$. For $1\le j \le N$, let $\vec f_j$
denote the corresponding vector in $\C^N$, and let $F$ denote the $N\times N$ matrix whose $j$th column is
$\vec f_j$.  Then
\begin{equation}\label{permint}
\int_\sn \prod_{j=1}^N(f_j\circ \pi_j){\rm d}\mu =  {1\over N!}{\rm perm}(F)\ .
\end{equation}

Let $\|\cdot\|_p$ denote the $L^p$ norm on $(\sn,\mu)$, and note that
$$|\vec f_j| = \sqrt{N}\|f_j\circ \pi_j\|_2$$
so that (\ref{main}) is equivalent to
\begin{equation}\label{intform}
\int_\sn \prod_{j=1}^N(f_j\circ \pi_j){\rm d}\mu  \le \prod_{j=1}^N\|f_j\circ \pi_j\|_2\ .
\end{equation}
In the from (\ref{intform}), the inequality (\ref{main}) bears a striking resemblance to the Brascamp--Lieb type  inequality  on
$S^{N-1}$ that we proved in proved in \cite{CLL}. For purposes of comparison, we recall this result.

Let $\nu$ denote the uniform probability measure on $S^{N-1}$, the unit sphere in $\R^N$. 
For each $j=1,2,\dots,N$, let $\vec e_j$ denote the $j$th standard basis vector in $\R^N$, so that for any
$\vec x$ in $S^{N-1}$, $\pi_j(\vec x) = \vec x\cdot \vec e_j$ is the $j$th component of $\vec x$. Then:

\bigskip
\begin{thm} \label{thm2}
For all $N\ge 2$, given non--negative  measurable functions
$f_1,\dots,f_N$, on $[-1,1]$, 
\begin{equation}\label{bound}
\int_{S^{N-1}}\left(\prod_{j=1}^N f_j\circ \pi_j\right){\rm d}\nu \le
\prod_{j=1}^N \|f_j\circ \pi_j\|_{L^p(S^{N-1})}\ .
\end{equation}
for all $p \geq 2$. Moreover, the $L^2$ norm  is optimal in that for each $p<2$,
there exist functions $f_j$  so that 
$\|f_j\circ\pi_j\|_{L^p(S^{N-1})} < \infty$ for each $j$, while the integral on the left side of
(\ref{bound}) diverges. Finally, for every $p \geq 2$ and $N\ge 3$, there is equality in (\ref{bound}) if and only if 
some function $f_j$ vanishes identically, or else each $f_j$ is constant.
\end{thm}
\bigskip

Note that Theorem \ref{thm2} provides  sharp information on the ratio
\begin{equation}\label{sphereratio}
\frac{\int_{S^{N-1}}\left(\prod_{j=1}^N f_j\circ \pi_j\right){\rm d}\nu}
{\prod_{j=1}^N \|f_j\circ \pi_j\|_{L^p(S^{N-1})}}
\end{equation}
for {\em all} values of $p$. However, Theorem \ref{thm1} only  provides sharp information on the ratio
\begin{equation}\label{symratio}
\frac{\int_\sn \prod_{j=1}^N(f_j\circ \pi_j){\rm d}\mu}
{\prod_{j=1}^N\|f_j\circ \pi_j\|_p}
\end{equation}
for $p\ge 2$.  Moreover, while the ratio in (\ref{sphereratio}) can have an infinite numerator, and finite denominator for $p<2$, this is clearly impossible for (\ref{symratio}). In fact, it is easy to obtain a sharp
upper bound on   (\ref{symratio}) for $p=1$.   

The fact that $2$ is the criitical $L^p$ index in the inequality (\ref{bound}) for all values of $N$ has interesting consequences for
the subadditivity of entropy on the sphere $S^N$, as explained in \cite{CLL}, and analogous entropy inequalities for
$\sn$ would follow from (\ref{intform}) in the same way. A recent paper \cite{BCM} of Barthe, Cordero-Erausquin and Maurey provides an illuminating geometric perspective on the  criticality of $p=2$ in  (\ref{bound}), among other things.

Thus, while (\ref{intform}) bears a clear family resemblance to (\ref{bound}), there
are differences.  Nonetheless, (\ref{intform}) can be proved by the same means that
were employed in \cite{CLL} to prove (\ref{bound}), as we explain in Section 2.
Next, in Section 3, we give an alternative proof of Theorem \ref{thm1}. In Section 4, we prove
the bounds on (\ref{symratio}) for $1\le p \le 2$, and we discuss certain 
natural conjectures and open problems. Finally, in an appendix we prove an interpolation theorem
that is used in Section 4.

%%%%%%%%%%%%%%%%%%%%%%%%%%%%%%%%%%%%%%%%%%%%%%%%%%%%%%%%%%%%%%%%%%%%%%%
\section{First Proof of Theorem \ref{thm1}}\label{fp}

Theorem \ref{thm2} was proved using a monotone heat kernel interpolation argument. Theorem \ref{thm1} can be proved in the same manner.

For $1\le i,j\le N$ with $i\ne j$ define the operator $D_{i,j}$ on $L^2(\sn,\mu)$ by
$$D_{i,j}g(\sigma) =   g(\sigma\sigma_{i,j}) - g(\sigma)\ .$$
One easily sees that each $D_{i,j}$ is  self  adjoint, and indeed, that
\begin{equation}\label{square2}
D_{i,j}^2 = -2D_{i,j}\ ,
\end{equation}
so that $D_{i,j}$ is actually non positive. The Laplacean $\Delta$ on $\sn$ is the operator
$$-\Delta = \sum_{i<j}D_{i,j}^2\ .$$
Note that a function $g$ on $\sn$ is of the from $f\circ \pi_j$ for some $f:\{1,\dots,N\}\to \R$ if and only if
$$D_{k,\ell} g= 0$$
whenever $k\ne j$ and $\ell \ne j$. 

Next observe that the left hand side of (\ref{main}) is only increased if we replace each entry in each $\vec f_j$ by its absolute value, and the right hand side is unchanged by this operation. Hence, in proving the inequality, we may assume without loss of generality that each of the functions $f_j$ is non negative, and of course, that none are identically zero.

For any non negative function $f$ on  $\{1,\dots,N\}$ and any $1\le j\le N$, and any $1 \le p < \infty$, consider the function  defined by
$ \left(e^{t\Delta} (f\circ \pi_j)^p\right)^{1/p}$.
Since $\Delta$ commutes with every $D_{k,\ell}$, whenever  $k\ne j$ and $\ell \ne j$ we have that
$$D_{k,\ell}\left(e^{t\Delta}(f\circ \pi_j)^p\right) = e^{t\Delta} \left(D_{k,\ell}(f\circ \pi_j)^p\right) = 0\ .$$
Hence $\left(e^{t\Delta} (f\circ \pi_j)^p\right)^{1/p}$ depends on $\sigma$ only through $\pi_j(\sigma)$, and for $t\ge 0$, we can define the function $f(\cdot,t)$ on $\{1,\dots, N\}$ by 
\begin{equation}\label{fdef}
f(t,\pi_j) = \left(e^{t\Delta}(f\circ \pi_j)^p\right)^{1/p}\ .
\end{equation}

Notice that since  $e^{t\Delta}$ preserves the integrals of functions,
\begin{equation}\label{const}
\|f(t,\pi_j)\|_p = \|f(0,\pi_j)\|_p = \|f\circ \pi_j\|_p\ 
\end{equation}
for all $t>0$.  Moreover,  since the null space of $\Delta$ consists of the constant functions on $\sn$ and nothing else,
and since all non zero eigenvalues of $\Delta$ are strictly negative, 
$\lim_{t\to \infty} e^{t\Delta}(f\circ \pi_j)^p = \int_\sn  (f\circ \pi_j)^p{\rm d}\mu$, and hence
\begin{equation}\label{conlim}
\lim_{t\to\infty}f(t,\pi_j) =  \|f\circ \pi_j\|_p\ .
\end{equation}

Now, given $N$ non negative, non zero  functions $f_j$ on $\{1,\dots,N\}$,  and $t>0$, define
\begin{equation}\label{etadef}
\eta_p(t) = \int_\sn \left(\prod_{j=1}^N f_j(t,\pi_j)\right){\rm d}\mu\ .
\end{equation}
(Note that by (\ref{fdef}), the right hand side does depend on $p$, as indicated by the subscript on the left hand side.)
By (\ref{conlim}), $\lim_{t\to\infty}\eta_p(t) = \prod_{j=1}^N\|f_j\circ \pi_j\|_p$, while clearly
$\eta_p(0) = \int_\sn\prod_{j=1}^Nf_j\circ \pi_j{\rm d}\mu$.  Hence the inequality
(\ref{intform}) would be proved if we could show that $t\mapsto \eta_2(t)$ is non decreasing.
Moreover, it is clear that for all $s,t>0$ and all $j$,
$$f(t+s,\pi_j) =  \left(e^{s\Delta} (f(t, \pi_j))^p\right)^{1/p}\ .$$
Because of this semigroup property, it would suffice to show that 
\begin{equation}\label{key}
{{\rm d}\over {\rm d}t}\eta_2(t)\bigg|_{t=0}\ge 0\ .
\end{equation}
This is indeed what we shall do. (Notice that the differentiability of $\eta_2$ is not an issue in this finite dimensional setting.)  The following lemmas prepare the way for the computation of the left side of (\ref{key}).

\medskip
\begin{lm}\label{lem1}  For any function $g$ on $\sn$,
\begin{equation}\label{leibniz}
\Delta g^2 = 2(\Delta g)g + 2\sum_{i<j}\left|D_{i,j}g\right|^2\ .
\end{equation}
\end{lm}
\medskip
\noindent{\bf Proof:} We compute
$$\left|D_{i,j}g(\sigma)\right|^2 = g^2(\sigma\sigma_{i,j}) + g^2(\sigma) -2g(\sigma\sigma_{i,j})g(\sigma)$$
and
$$(2D_{i,j}g(\sigma))g(\sigma) = 2g(\sigma\sigma_{i,j})g(\sigma) - 2g^2(\sigma)\ .$$
Summing these equations, one has
\begin{eqnarray}
\left|D_{i,j}g(\sigma)\right|^2 + (2D_{i,j}g(\sigma))g(\sigma) &=& 
 g^2(\sigma\sigma_{i,j}) - g^2(\sigma)\nonumber\\
 &=& D_{i,j}g^2(\sigma)\ .\nonumber\\
 \end{eqnarray}
 Multiplying through by 2, summing over $i<j$, and making use of (\ref{square2}), one obtains (\ref{leibniz}). \qed
 \medskip
 
 Lemma \ref{lem1} justifies the following notational convention: We define $|\nabla g|^2$
 by
 \begin{equation}\label{nablasq}
 |\nabla g|^2 =  \sum_{i<j}\left|D_{i,j}g\right|^2\ .
 \end{equation}

 \medskip
\begin{lm}\label{lem2}
With $f(t,\sigma)$  defined by 
$f(t,\sigma) = \left(e^{t\Delta} (f\circ \pi_j)^2(\sigma)\right)^{1/2}$,
 \begin{equation}\label{deriv}{\partial \over \partial t}f(t,\sigma)\bigg|_{t=0} = \Delta (f\circ \pi_j) + {\left|\nabla(f\circ \pi_j)\right|^2 \over
(f\circ \pi_j)}\ .
\end{equation}
\end{lm}
\medskip
\noindent{\bf Proof:} This is a simple computation using (\ref{leibniz}). \qed
\medskip

The following lemma  gives the
modification to the Leibniz rule for the finite difference operation $D_{i,j}$. What is crucial for us is that the modification drops out if one of the functions does not depend on $\sigma$ through either $\pi_i(\sigma)$ or 
$\pi_j(\sigma)$.

\medskip
\begin{lm}\label{lem3}
{\it   For any two function $g$ and $h$ on $\sn$, and any $1\le i<j\le N$,
$$D_{i,j}(gh)(\sigma) = \left(D_{i,j}g(\sigma)\right)h(\sigma) + g(\sigma\sigma_{i,j})\left(D_{i,j}h(\sigma)\right)\ .$$
In particular, if $D_{i,j}h = 0$,
 \begin{equation}\label{leibniz2}
 D_{i,j}(gh) = \left(D_{i,j}g\right)h\ .
 \end{equation}}
 \end{lm}
\medskip
\noindent{\bf Proof:} This is an even  simpler computation, in which one makes the  
obvious addition and subtraction. \qed

\medskip
\begin{lm}\label{lem4}
The inequality (\ref{intform}) is satisfied for any $N$ non negative functions $f_j$ on $\{1,\dots,N\}$
 \end{lm}
\medskip
\noindent{\bf Proof:} We may freely assume that none of the functions is identically zero, since then 
(\ref{intform}) is trivially satisfied, with zero on both sides.

Define $f_j(t,\cdot)$ by
$$f_j(t,\pi_j) = \left(e^{t\Delta} (f_j\circ \pi_j)^2\right)^{1/2}\  .$$
and $\eta_2(t)$ by
$\eta_2(t) = \int_\sn \left(\prod_{j=1}^N f_j(t,\pi_j)\right){\rm d}\mu$.
Then by Lemma \ref{lem2},
$${{\rm d}\over {\rm d}t}\eta_2(t)\bigg|_{t=0} =
\sum_{j=1}^N \int_\sn \left(\Delta (f_j\circ \pi_j) + {|\nabla (f_j\circ \pi_j)|^2\over
(f_j\circ \pi_j)}\right)  \left(\prod_{k=1,k\ne j}^N (f_k\circ \pi_k)\right){\rm d}\mu\ .$$
Consider the contribution coming from
$$\int_\sn \Delta (f_j\circ \pi_j) \left(\prod_{k=1,k\ne j}^N (f_k\circ \pi_k)\right){\rm d}\mu\ .$$
Notice that
$$-\Delta (f_j\circ \pi_j)   = \sum_{i>j}D_{j,i}^2(f_j\circ \pi_j) +  \sum_{i< j}D_{i,j}^2(f_j\circ \pi_j)\ .$$
Since $D_{i,j}$ is self adjoint,
\begin{eqnarray}
-\int_\sn \Delta (f_j\circ \pi_j) \left(\prod_{k=1,k\ne j}^N (f_k\circ \pi_k)\right){\rm d}\mu &=&
\sum_{i>j}\int_\sn D_{j,i}(f_j\circ \pi_j) D_{j,i}\left(\prod_{k=1,k\ne j}^N (f_k\circ \pi_k)\right){\rm d}\mu\nonumber\\
&+& \sum_{i<j}\int_\sn D_{i,j}(f_j\circ \pi_j) D_{i,j}\left(\prod_{k=1,k\ne j}^N (f_k\circ \pi_k)\right){\rm d}\mu\ .\nonumber\\
\end{eqnarray}
Now by (\ref{leibniz2}),
\begin{eqnarray}
&\phantom{=}&\int_\sn D_{j,i}(f_j\circ \pi_j) D_{j,i}\left(\prod_{k=1,k\ne j}^N (f_k\circ \pi_k)\right){\rm d}\mu\nonumber\\
&=& \int_\sn D_{j,i}(f_j\circ \pi_j) D_{j,i}(f_i\circ\pi_i)\left(\prod_{k=1,k\ne j,i}^N (f_k\circ \pi_k)\right){\rm d}\mu\nonumber\\
&=& \int_\sn {D_{j,i}(f_j\circ \pi_j)\over f_j\circ \pi_j}{ D_{j,i}(f_i\circ\pi_i)\over f_i\circ \pi_i} 
\left(\prod_{k =1}^N (f_k\circ \pi_k)\right){\rm d}\mu\ .\nonumber\\
\end{eqnarray}

Defining the non negative function $\rho$  by
$\rho = \prod_{k =1}^N (f_k\circ \pi_k)$
and making a similar computation for the sum on $i<j$, we obtain
\begin{eqnarray}
&\phantom{=}&\int_\sn \Delta (f_j\circ \pi_j) \left(\prod_{k=1,k\ne j}^N (f_k\circ \pi_k)\right){\rm d}\mu=\nonumber\\
&\phantom{=}&\sum_{i\ne j}\int_\sn\left( {D_{i,j}(f_j\circ \pi_j)\over f_j\circ \pi_j}{ D_{i,j}(f_i\circ\pi_i)\over f_i\circ \pi_i}\right) \rho{\rm d}\mu\ .\nonumber\\
\end{eqnarray}

{}From here we see that
\begin{equation}\label{square}
{{\rm d}\over {\rm d}t}\eta(t)\bigg|_{t=0} =
{1\over 2} \sum_{i,j} \int_\sn \left(
 {D_{i,j}(f_j\circ \pi_j)\over f_j\circ \pi_j} - { D_{i,j}(f_i\circ\pi_i)\over f_i\circ \pi_i} \right)^2 \rho{\rm d}\mu\ .
 \end{equation}    \qed
 
\medskip
\begin{lm}\label{lem5} If there is equality in (\ref{intform}), and if none of the functions $f_j$ is identically zero, then
each of them has a constant modulus. That is, for each $j$ and $k$, $|f_j(k)| = \|f_j\|_2$.
 \end{lm}
\medskip
\noindent{\bf Proof:}   Suppose that for some $N$ functions $f_j$ on $\{1,\dots,N\}$, there is equality in 
(\ref{intform}), and that none of the functions vanishes identically.  Then clearly there is still equality in
in  (\ref{intform}) if we replace each $f_j$ by $|f_j|$. Hence, we may freely assume that the functions are all non negative, 
and that none vanishes identically. 

It now follows from the proof of Lemma \ref{lem4} that
${\displaystyle{{\rm d}\over {\rm d}t}\eta_2(t) = 0}$ for all $t$.  However, for all $t>0$,
each $f_j(t,\cdot)$ will be strictly positive, and so from (\ref{square}) we see that for each $t>0$,
we must have
\begin{equation}\label{cases} 
{D_{i,j}(f_j(t, \pi_j)\over f_j(t, \pi_j)} = { D_{i,j}(f_i(t,\pi_i)\over f_i(t, \pi_i)} 
\end{equation}
for all $i$ and $j$, at every $\sigma$. 

Fix any $i\ne j$ in $\{1,\dots,N\}$. For $N>2$, we can chose $k$ from $\{1,\dots,N\}$
so that $k \ne i,j$. Note that
$$\pi_i(\sigma_{i,k}\sigma_{i,j}) = \pi_i(\sigma_{i,j}) = j\ .$$
Thus, the right side of (\ref{cases}) vanishes at $\sigma = \sigma_{i,k}$.
Hence the left side vanishes there as well, and since 
$$\pi_j(\sigma_{i,k}\sigma_{i,j}) = k\qquad{\rm while}\qquad  \pi_j(\sigma_{i,j}) = i\ .$$
We therefore conclude that for each $j$, $f_j(t,k) = f_j(t,i)$ for all $i,k\ne j$. 

This is almost what we seek. To conclude the proof, consider the permutation $\sigma_c$
with
$$\sigma_c(i) = j\qquad  \sigma_c(j) = k \qquad{\rm and}\qquad \sigma_c(k) = i\ .$$
Since, 
$$\pi_j(\sigma_{c}\sigma_{i,j}) = k\qquad{\rm and}\qquad  \pi_j(\sigma_{c}) = j\ ,$$
the numerator in the right side of (\ref{cases})  is $f_i(t,k) - f_i(t,j)$ at $\sigma_c$. 
By what we have see above, this is zero. Therefore, the numerator of the left hand side of (\ref{cases})
vanishes at $\sigma_c$. This is $f_j(t,j) - f_j(t,i)$. Hence for each $t>0$, and each $j$, 
$f_j(t,\cdot)$ is constant. By  continuity, it follows that each $f_j(\cdot)$ is constant. \qed

To complete the proof of Theorem \ref{thm1}, consider $N$ functions $f_j$ on $\{1,\dots,N\}$, none of which is identically zero, and for which equality holds (\ref{intform}). Then we know that each $|f_j(k)|$ is non zero, and so we can define a complex number $Z_{k,j}$ by
$$Z_{k,j} = f_j(k)/|f_j(k)|\ .$$
Clearly, each of the $Z_{j,k}$ lies on the unit circle. 

Now let $A$ be the $N\times N$ matrix with $A_{j,k} = |f_{j,k}|$, and let $Z$ be the $N\times N$
matrix with entries $Z_{j,k} $. Then, if $F$ is the $N\times N$ matrix with $F_{j,k} = f_j(k)$,
$F = Z\cdot A$, where the write hand side is the Hadamard product of $Z$ and $A$. 
Clearly, 
\begin{equation}
|{\rm perm}(Z\cdot A)| \le  |{\rm perm}(A)| 
\end{equation}
and there is equality if and only if the quantity
\begin{equation}\label{uni}
\prod_{j=1}^N Z_{j,\sigma(j)}
\end{equation}
does not depend on $\sigma$.

\medskip
\begin{lm}\label{lem6}  Let $Z$ be an $N\times N$ matrix such 
that for each $j$ and $k$, the $j,k$th entry is
  a complex number $Z_{k,j}$ lying on the unit circle.  Then, the
  product in (\ref{uni}) is independent of $\sigma$ if and only if
  there are vectors $\vec \xi$ and $\vec \zeta$ in $\C^N$ with each
  entry lying in the unit circle such that for each $j$ and $k$,
  $Z_{j,k} = \xi_j\zeta_k$.
 \end{lm}
\medskip
\noindent{\bf Proof:}   
Suppose that   $Z_{j,k} = \xi_j\zeta_k$.  Then 
$$\prod_{j=1}^N Z_{j,\sigma(j)} = \prod_{j=1}^N \xi_j\zeta_{\sigma(j)} = \left(\prod_{j=1}^N\xi_j\right)
 \left(\prod_{k=1}^N\zeta_k\right)\ ,$$
 which is independent of $\sigma$. This proves sufficiency.
 
To prove necessity, we use induction. The lemma is clearly
true for $N\leq 2$, so we start with 3.   Let us expand in the first row of
$Z$.
$$
{\rm perm}(Z) = \sum_{i=1}^N Z_{1i} \, {\rm perm}(Z^i)    
$$
where $Z^i$ is the the $(N-1)\times (N-1)$ matrix cofactor of $Z_{1i}$ in $Z$.  
Since $|{\rm perm}(Z)| =N!$ and
$|{\rm perm}(Z^i )|\leq (N-1)!$ we must have
$$
|{\rm perm}(Z^i)| = (N-1)! \ .
$$
By induction,  $Z^i$ must be of the form
$$
Z_{j,k} = \xi_j\zeta_k   \qquad  {\rm for}\qquad  j \neq 1,  \  k \neq i\ . 
$$

Likewise any $(N-1)\times (N-1)$ submatrix complementary to $Z_{\alpha, \beta}$
must have the form
\begin{equation}  \label{zee}
Z_{j,k} =  \xi_j\zeta_k   \qquad  {\rm for} j\neq \alpha,
\quad k\neq \beta\ . 
\end{equation}

The $(N-1)$ dimensional vectors $\xi$ and $\zeta$ depend on $\alpha $ and
$\beta$, in principle, but this is not so, as we now show. Let $1\leq j,\ k, \ l,\ m \leq N$ be four 
integers. There is an $\alpha \leq N$ that is different from $j$ and $l$. Similarly, there is a 
$\beta$ different from $k$ and $m$. Then, equation (\ref{zee}) is valid,  with the {\it same}
$\xi$ and $\zeta$,   for both
$Z_{j,k} $ and for $Z_{l,m} $, that is,  $Z_{j,k} = \xi_j\zeta_k $ and $Z_{l,m} = \xi_l\zeta_m$.
From this we see that  
$$
Z_{j,k}  Z_{l,m }  =   Z_{j,m}  Z_{l,k}  \ 
$$
for {\it any} quartet of indices $1\leq j,\ k, \ l,\ m \leq N$.  With $l=m=1$ we then deduce 
that   $ Z_{j,k} = Z_{j,1} Z_{1,k}  /  Z_{1,1}$, and we are done. 
\qed
\medskip

\noindent{\bf Proof of Theorem \ref{thm1}}  We have already explained 
that (\ref{intform}) is equivalent to (\ref{main}),
and (\ref{intform}) has been proved in Lemma \ref{lem4}. The statement
concerning the cases of equality then follows from Lemmas \ref{lem5}
and \ref{lem6}. \qed
\medskip
%%%%%%%%%%%%%%%%%%%%%%%%%%%%%%%%%%%%%%%%
%%%%%%%%%%%%%%%%%%%%%%%%%%%%%%%%%%%%%%%%
%%%%% BEGINNING OF SECTION 3

\section{Second Proof of Theorem \ref{thm1}}\label{sp}

\bigskip
The second proof is based on induction and the arithmetic--geometric mean inequality.
In that sense, it uses only elementary tools. However, as will be seen, they must be applied in a particularly judicious way. In any case, 
the structure of this proof leads naturally to a generalization of Theorem \ref{thm1}
to non square matrices $F$. 

Given $K \le N$ vectors in $\C^N$, we from the $K\times N$ matrix whose $i$th  {\em row}
is $\vec f_i$. Let $f_{i,j}$ denote the $i,j$th entry of this matrix; i.e., the $j$th
 entry of $\vec f_i$. Define the functional ${\cal P}(\vec f_1, \dots, \vec f_K)$ by
\begin{equation}\label{pdef}
{\cal P}(\vec f_1, \dots, \vec f_K) =\left[ \sum_{1 \le j_1 < j_2 < \dots <j_K \le N}
\left({\rm perm} \begin{bmatrix} f_{1,j_1} & f_{1,j_2} & \dots & f_{1,j_K} \\ 
                      f_{2,j_1} & f_{2,j_2} & \dots & f_{2,j_K} \\
                      \dots & \dots & \dots  & \dots \\
                      f_{K,j_1} & f_{K,j_2} & \dots & f_{K,j_K}
                      \end{bmatrix} \right)^2\right]^{1/2} \ .
                      \end{equation}
                      
Notice that each of the permanents in (\ref{pdef}) is the permanent of a $K\times K$ matrix. Note also that
for $K = N$, there is only one term in the sum, and
\begin{equation}\label{peq}
{\cal P}(\vec f_1, \dots, \vec f_N)  = {\rm perm}[\vec f_1, \dots, \vec f_N]\ .
\end{equation}
We shall now prove:

\bigskip
\begin{thm}\label{thm4} The inequality 
\begin{equation}\label{pineq}
{\cal P}(\vec f_1, \dots, \vec f_K) \le \sqrt{{N \choose K}} {K!\over N^{K/2}}\prod_{j=1}^K|\vec f_j|\ .
\end{equation}
If $K \ge 2$ and none of the vectors $\vec f_i$ is the zero vector, then (\ref{pineq}) holds with equality if and
only if $[\vec f_1, \dots, \vec f_K]$ is a rank one matrix, and each the vectors $\vec f_i$ is a constant modulus vector. 
\end{thm}
\medskip

Because of (\ref{peq}) Theorem \ref{thm4} reduces to Theorem \ref{thm1}  in the case $K=N$.

\medskip
\noindent{\bf Proof:} As before, when proving the inequality, we may assume that all entries of each vector are non negative.
The proof proceeds by induction in $K$. The inequality is trivial in case $K=1$,
although any vector yields equality in this case. The first non trivial case is $K=2$. We now treat this case explicitly, since the same sort of reasoning will be employed in the general inductive step.

\begin{equation}\label{p21}
{\cal P}( \vec f_1, \vec f_2)^2 =\sum_{i < j} (f_{1,i}f_{2,j} + f_{1,j}f_{2,i})^2
=\sum_{i < j} f_{1,i}^2f_{2,j}^2 + f_{1,j}^2f_{2,i}^2 + \sum_{i < j} 2  f_{1,i}f_{2,j}f_{1,j}f_{2,i} \ .
\end{equation}
We will use the arithmetic--geometric mean inequality on the terms of the last sum, but there are two natural ways to do this, {\em and we need to use some of each}. Therefore, rewrite the last sum as
$$
\sum_{i < j} 2  f_{1,i}f_{2,j}f_{1,j}f_{2,i} = \alpha \sum_{i < j} 2 ( f_{1,i}f_{2,j})(f_{1,j}f_{2,i})
+(1-\alpha)  \sum_{i < j} 2 ( f_{1,i}f_{2,i})(f_{1,j}f_{2,j})
$$
where $\alpha <1$ will be determined later. 
By  the arithmetic--geometric mean  mean inequality,
\begin{equation}\label{p22}
\sum_{i < j} 2  f_{1,i}f_{2,j}f_{1,j}f_{2,i} \le \alpha \sum_{i < j} \left[ ( f_{1,i}f_{2,j})^2+(f_{1,j}f_{2,i})^2\right]
+(1-\alpha)  \sum_{i < j}  \left[( f_{1,i}f_{2,i})^2 + (f_{1,j}f_{2,j})^2\right] \ .
\end{equation}
Combing (\ref{p21}) and (\ref{p22}),
$$
{\cal P}( \vec f_1, \vec f_2)^2 \le (1 +\alpha)    \sum_{i \not= j} f_{1,i}^2f_{2,j}^2  +(1-\alpha)(N-1)  \sum_{i }  ( f_{1,i}f_{2,i})^2 \ .
$$
Choose $\alpha= (N-2)/N$ yields $1+\alpha = (1-\alpha)(N-1)=2(N-1)/N$,  and  find
$$
{\cal P}( \vec f_1, \vec f_2)^2 \leq  2 {N-1 \over N}  \sum_{i , j} f_{1,i}^2f_{2,j}^2 
= 2{(N-1) \over N} |\vec f_1|^2 |\vec f_2|^2 \ ,
$$
which is the desired result for $K=2$.

\medskip
For the general case we can write
$$
{\cal P}(\vec f_1, \dots, \vec f_K) ^2
= {1 \over K!} \sum'_{j_1, \dots , j_K} \left(\sum_{k=1}^K f_{1,j_k} a_{j_1, \dots , \widehat{j_k}, \dots j_K}\right)^2
$$
where the prime over the summation sign indicates that the sum is over distinct indices. Here
$$
a_{j_1, \dots , \widehat{j_k}, \dots j_K} ={\rm perm}
\begin{bmatrix}
f_{2,j_1} & \dots & \widehat{ f_{2,j_k}} &
\dots & f_{2,j_K} \\ \dots \\  
f_{K,j_1} & \dots & \widehat{ f_{K,j_k}} &
\dots & f_{K,j_K}
\end{bmatrix}
$$
is the permanent of the matrix where the
$f_1$ row and the $j_k$ column has been removed so that we take the permanent of an
$K-1 \times K-1$ matrix. The removal of the column is indicated by the circumflex symbol, as usual.
 Note that this expression is symmetric in the indices. Developing the square leads to
$$
 {1 \over K!} \sum'_{j_1, \dots , j_K} \sum_{k, l}^K f_{1,j_k} f_{1,j_l}a_{j_1, \dots , \widehat{j_k}, \dots j_K} a_{j_1, \dots , \widehat{j_l}, \dots j_K}
= {1 \over K!} \sum'_{j_1, \dots , j_K} \sum_{k=1}^K (f_{1,j_k})^2 (a_{j_1, \dots , \widehat{j_k}, \dots j_K} )^2 
$$
$$
+  {1 \over K!} \sum'_{j_1, \dots , j_K} \sum_{k \not= l}^K  f_{1,j_k} f_{1,j_l}a_{j_1, \dots , \widehat{j_k}, \dots j_K} a_{j_1, \dots , \widehat{j_l}, \dots j_K}
$$
As before, we estimate the last term using the  the arithmetic--geometric mean inequality {\em in two different ways} and
obtain
\begin{equation}\label{p25}
\alpha  {1 \over (K-2)!} \sum'_{j_1, \dots , j_K}  ( f_{1,j_1})^2 (a_{j_2, \dots j_K})^2 
+(1-\alpha) {1 \over (K-2)!} \sum'_{j_1, \dots , j_K}  (f_{1,j_1})^2 (a_{j_1, \widehat{j_2}, \dots j_K})^2 \ 
\end{equation}
where $0 < \alpha < 1$ will be chosen below.  

Focusing on the very last term, we  write it as
\begin{eqnarray}
&\phantom{=}&(1-\alpha)  {1 \over (K-2)!} \sum'_{j_1, \dots , j_K}   (f_{1,j_1})^2 (a_{j_1,\widehat{j_2}, \dots j_K})^2\nonumber\\
&=&(1-\alpha)(N-K+1)  {1 \over (K-2)!} \sum_{j_1} (f_{1,j_1})^2  \sum'_{j_2, \dots , j_{(K-1)} \not= j_1}   (a_{j_1, \dots j_{(K-1)}})^2
\nonumber\\
&=&(1-\alpha){(N-K+1) \over K-1}  {1 \over (K-2)!} \sum_{j_1} (f_{1,j_1})^2 
\sum_{l=1}^{K-1} \sum'_{j_2, \dots , j_{(K-1)} \not= j_1}   (a_{j_2, \dots j_1, \dots j_{(K-1)}})^2
\nonumber\\
\end{eqnarray}
where the index $j_1$ is in the $l$-th position.  Note that
$$\sum_{l=1}^{K-1} \sum'_{j_2, \dots , j_{(K-1)} \not= j_1}   (a_{j_2, \dots j_1, \dots j_{(K-1)}})^2
= \sum'_{j_2, \dots , j_K \not= j_1} (a_{j_2, \dots j_K} )^2 \ .$$
Hence, collecting the terms we get
\begin{eqnarray}
{\cal P}(\vec f_1, \dots, \vec f_K) ^2
&\leq& (1+(K-1) \alpha)   {1 \over (K-1)!} \sum_{j_1} ( f_{1,j_1})^2 \sum'_{j_2, \dots , j_K \not= j_1} (a_{j_2, \dots j_K} )^2 \nonumber\\
&+& (1-\alpha){(N-K+1)}  {1 \over (K-1)!} \sum_{j_1}( f_{1,j_1})^2 
\sum'_{j_2, \dots , j_K \not= j_1} (a_{j_2, \dots j_K} )^2 \nonumber\\
\end{eqnarray}

Now choose $\alpha = (N-K)/N$ so that
$$
(1+(K-1) \alpha)=(1-\alpha){(N-K+1)} = {K(N-K+1) \over N}\ ,
$$
and
$$
{\cal P}(\vec f_1, \dots, \vec f_K) ^2
\le {K(N-K+1) \over N}{1 \over (K-1)!} 
\sum_{j_1} ( f_{1,j_1})^2 \sum'_{j_2, \dots , j_K} (a_{j_2, \dots j_K} )^2 
$$
By the  inductive hypotheses,
$$
\sum'_{j_2, \dots , j_K} (a_{j_2, \dots j_K} )^2 \le (K-1)! {N \choose K-1} {((K-1)!)^2 \over N^{K-1}}
\Pi_{j=2}^K |\vec f_j| ^2 \ ,
$$
which yields
$$
{\cal P}(\vec f_1, \dots, \vec f_K) ^2
\le { N \choose K} {(K!)^2 \over N^K} \Pi_{j=1}^N |\vec f_j |^2 \ .
$$
This proves the inequality.

To establish the cases of equality, we return to the step where the arithmetic--geometric mean
inequality has been employed, of course still under the assumption that the entries are all non negative. Equality entails that
$$
f_{1,j}f_{2,j} = f_{1,i}f_{2,i} \ {\rm and} \ f_{1,i}f_{2,j} = f_{1,j}f_{2,i} \ , 
$$
for all $i\ne j$.
Since the vectors $\vec f_1$ and $\vec f_2$ are not allowed to be the zero vectors there exists
an index $i$ so that $f_{1,i}$ is not zero. If $f_{2,i} = 0$ then it follows from $f_{1,i}f_{2,j} = f_{1,j}f_{2,i}$ that
$f_{2,j}=0$ for all $j \not= i$. This would mean that $f_2$ is identically zero, contrary to the conditions
in the theorem. Thus, both, $f_{1,i}$ and $f_{2,i}$
are non zero. It now follows from  $f_{1,j}f_{2,j} = f_{1,i}f_{2,i}$  that neither $f_{1,j}$ nor $f_{2,j}$ can
be zero. Therefore
$$
{f_{2,j} \over f_{2,i}} = {f_{1,j} \over f_{1,i}} = {f_{2,i} \over f_{2,j}}
$$
for all $i < j$ and hence the vectors $\vec f_1$ and $\vec f_2$ have to be constant
vectors. Since the same argument applies to any two vectors the result follows. 
We may now treat the case that the entries are complex using the final lemma of the previous section.
\qed

\bigskip
\begin{cl}\label{cor1} 
For all $p$ with $1\le p \le 2$,
the quantity
\begin{equation}
{\cal P}_p(\vec f_1, \dots, \vec f_K) =\left[ \sum_{1 \le j_1 < j_2 < \dots <j_K \le n}
\left({\rm perm} 
\begin{bmatrix}\label{56}
 f_{1,j_1} & f_{1,j_2} & \dots & f_{1,j_K} \\ 
                      f_{2,j_1} & f_{2,j_2} & \dots & f_{2,j_K} \\
                      \dots & \dots & \dots  & \dots \\
                      f_{k,j_1} & f_{k,j_2} & \dots & f_{K,j_K} 
                      \end{bmatrix}\right)^p
                      \right]^{1/p} 
\end{equation}
satisfies the inequality
$$
{\cal P}_p(\vec f_1, \dots, \vec f_K) \le {N \choose K}^{1/p}  {K!\over N^{K/2}}\prod_{j=1}^K|\vec f_j| \ .
$$
\end{cl}
\medskip

\noindent{\bf Proof:} 
By  H\"older's inequality,
$$
{\cal P}_p(\vec f_1, \dots, \vec f_K) \le  {N \choose K}^{1/p-1/2} {\cal P}(\vec f_1, \dots, \vec f_K)\ .
$$
Now apply Theorem \ref{thm4} to estimate ${\cal P}(\vec f_1, \dots, \vec f_K)$; this yields the result.
\qed

\medskip
%%%%%%%%%%%%%%%%%%%%%%%%%%%%%%%%%%%%%%%%
%%%%%%%%%%%%%%%%%%%%%%%%%%%%%%%%%%%%%%%%
%%%%% BEGINNING OF SECTION 4

\section{Bound for other values of $p$}\label{otherp}

For $1\le p < \infty$, and any vector $\vec f$ in $\C^N$, define
\begin{equation}\label{plength}
|\vec f|_p = \left(
\sum_{k=1}^N |(\vec f_j)_k|^p\right)^{1/p} \ .
\end{equation}
Note that if $f$ is the corresponding function of $\{1,\dots,N\}$, for each $j$ we have
\begin{equation}\label{prel}
|\vec f|_p =N^{1/p}\|f\circ \pi_j\|_p\ .
\end{equation}
By (\ref{permint}) and (\ref{prel}),
\begin{equation}\label{otherratio}
\frac{\int_\sn \prod_{j=1}^N(f_j\circ \pi_j){\rm d}\mu}
{\prod_{j=1}^N\|f_j\circ \pi_j\|_p} = 
\frac{N!}{N^{N/p}} \frac{|{\rm perm}[\vec f_1, \dots, \vec f_N]| }{ \prod_{j=1}^N | \vec f_j|_p} \ .
\end{equation}
Thus, we may study the ratio in (\ref{symratio}) by studying the ratio on the right in
(\ref{otherratio}).

Define the function
\begin{equation}\label{pratio}
C(p) =  
\sup_{\vec f_1,\dots,\vec f_N \ne 0} \left\{ {|{\rm perm}[\vec f_1, \dots, \vec f_N]| \over \prod_{j=1}^N | \vec f_j|_p}  \right\} \ . 
\end{equation}
We know from Theorem 1 that
\begin{equation}\label{2ratio}
C(2) ={N! \over N^{N/2}} 
\end{equation}
with equality precisely when $[\vec f_1, \dots, \vec f_N]$ is the constant matrix.
Moreover, it is easy to see that $C(1) =1$:
Observe that
\begin{equation}\label{57}
|{\rm perm}[\vec f_1, \dots, \vec f_N]| \le \sum_{\sigma \in \sn} \prod_{j=1}^N|f_{j, \sigma(j)}| 
\end{equation}
and that
\begin{equation}\label{58}
\prod_{j=1}^N | \vec f_j| = \sum_{k_1, \dots, k_N} \prod_{j=1}^N  |f_{j, k_j}| 
\end{equation}
and note that every term in (\ref{57}) shows up in (\ref{58}), and hence
$C(1) \le  1$.
Choosing $\vec f_j = \vec e_j$ for each $j$, so that $F$ equals the identity matrix, shows that 
\begin{equation}\label{59}
C(1)=1 \ . 
\end{equation} 
In fact, there is equality only if
in each row and each column  of $F$ there is a single non-zero entry.

Notice that the optimizers for $C(p)$ are different for $p=1$ and $p=2$: For $p=1$ we get the optimal ratio by taking $\vec f_j = \vec e_j$ for each $j$, while for $p=2$ we get the optimal ratio by taking $\vec f_j$ to be constant for each $j$. 

If $\vec f_j = \vec e_j$ for each $j$, then 
$$ {|{\rm perm}[\vec f_1, \dots, \vec f_N]| \over \prod_{j=1}^N | \vec f_j|_p}  = 1\ .$$

However, if $\vec f_j$ is a non zero constant vector for each $j$, then
$$ {|{\rm perm}[\vec f_1, \dots, \vec f_N]| \over \prod_{j=1}^N | \vec f_j|_p}  = {N!\over N^{N/p}}\ .$$

Evidently,
\begin{equation}\label{cap}
C(p) \ge \max\left\{\ 1\ , \ {N!\over N^{N/p}}\ \right\}\ .
\end{equation}
Note that there is equality at $p=1$ and $p=2$. Pietro Caputo, to whom we sent an early draft of this paper, has suggested to us that  in fact there should be equality in (\ref{cap})
for $1 < p < 2$ as well.  It is easy to see that this is true for $N=2$. In that case, we may assume without
loss of generality that $\vec f_1 = \begin{bmatrix}1\\ x\\ \end{bmatrix}$ and 
$\vec f_2 = \begin{bmatrix}y\\ 1\\ \end{bmatrix}$ for some non negative numbers $x$ and $y$.
Then ${\rm perm}[\vec f_1,\vec f_2] = 1 + xy$, and by H\"older's inequality,
$$1 +xy \le (1+x^p)^{1/p}(1+y^q)^{1/q}$$
where $1/q = 1 - 1/p$. For $1 \le p \le 2$, $q \ge p$, and so
$$(1+y^q)^{1/q} \le (1+y^p)^{1/p}\ ,$$
with equality for $p < 2 < q$ if and only if $y=0$. 
We conclude that
$$1 +xy \le (1+x^p)^{1/p}(1+y^p)^{1/p}\ ,$$
which is the desired inequality.
Moreover,  for $p=2$, by the condition for equality in the Schwarz inequality,
there is equality if and only if $x=y=1$, while when $1 \le p < 2$, there is equality if and only if $x=y=0$.

Beyond this trivial case,  have not succeeded in proving the conjecture, but we do have the following
upper bound on $C(p)$:

\medskip
\begin{thm}\label{thm5} The function $ln(C(p))$ is a  convex function of $1/p$.
In particualar, for all $1\le p \le 2$,
$$C(p) \le  \left( {N! \over N^{N/2}}\right)^{2-2/p}\ .$$
\end{thm}

\medskip  
\noindent{\bf Proof:} 
The  first statement follows directly from the interpolation theorem in the appendix. This is a version of the Riesz--Thorin \cite{T} interpolation theorem for multilinear forms. The usual proof of the Riesz--Thorin interpolation theorem
for operators is, in fact, an interpolation theorem for bilinear forms; see for example, \cite{G}. It easily extends to multilinear forms, and though this seems likely to be known, we have not found any reference. We therefore
include the short proof in the appendix.

Given the logarithmic convexity,  for 
$t$ solving the equation
$
t + (1-t)/2 = 1/p$,
$$
C(p) \le C(1)^t C(2)^{1-t}.
$$
The rest now follows from our computation of $C(1)$ and $C(2)$. \qed
\medskip

One might try to compute $C(p)$ for $1< p < 2$ by adapting either of the two computations 
we made for $p=2$. Unfortunately, we have not been able to adapt either one. For the second
computation, the trouble arises at the point where we ``develop the square''.  For the first, 
there is an even more fundamental problem: The interpolation used there simply {\em is not monotone} for $p < 2$. 

To see this consider  $N =3$ and the vectors 
$\vec f_1 = \begin{bmatrix}1\\ x\\  y \end{bmatrix}$, 
$\vec f_2 = \begin{bmatrix}y\\ 1\\  x\end{bmatrix}$ and $\vec f_3 = \begin{bmatrix}x\\ y\\ 1\end{bmatrix}$
for some numbers $x$ and $y$ with $0 \le x,y\le 1$.
Then
$$[\vec f_1,\vec f_2,\vec f_3] = \begin{bmatrix}
1 & y & x\\
x & 1 & y\\
y & x & 1\\
\end{bmatrix}\ ,$$
which is a circulant matrix.

Define the function $\phi(x,y)$ by
$$\phi(x,y) = \frac{{\rm perm}[\vec f_1,\vec f_2,\vec f_3]}{|\vec f_1|_p|\vec f_2|_p|\vec f_3|_p}
 = \frac{1 + x^3 + y^3 + 3xy}{(1+x^p+y^p)^{3/p}}\ .$$

It is easy to see that the class of circulant matrices  is preserved under the heat semigroup flow used in the proof of Theorem \ref{thm1} provided in Section 2, so that this flow corresponds to a path $(x(t),y(t))$
on the graph of $\phi$ over  the unit square $0 \le x,y\le 1$. This path starts from the given values of
$x$ and $y$ and satisfies
$$\lim_{t\to\infty}(x(t),y(t)) = (1,1)\ .$$

However, for $p<2$, $\phi(x,y)$ has a strict local maximum at both $(x,y) = (1,1)$ and $(x,y) = (0,0)$.
Thus, for $p<2$, the heat semigroup flow is initially strictly monotone decreasing when started from 
$\vec f_j = \vec e_j$, $j = 1,2,3$. 
Nonetheless, further analysis of the function $\phi$ supports the conjecture; the example simply shows  that no flow preserving the class of circulant matrices can be used to prove it.

%%%%%%%%%%%%%%%%%%%%%%%%%%%%%%%%%%%%%%%%
%%%%%%%%%%%%%%%%%%%%%%%%%%%%%%%%%%%%%%%%
%%%%% BEGINNING OF SECTION 4

\section{Appendix }\label{interp}

Here we prove the following multilinear generalization of the Riesz--Thorin interpolation theorem

Let $J$ denote a multilinear functional of $M$--tuples of vectors $\vec f_j$ in $\C^N$.
Then with $f_{j,k}$ denoting the $k$th component of the vector $\vec f_j$, there are numbers
$J_{k_1,\dots,k_M}$ such that
\begin{equation}\label{sumform}
J(\vec f_1,\dots,\vec f_M) = \sum_{k_1,\dots,k_M}
J_{k_1,\dots,k_M}\prod_{j=1}^M f_{j,k_j}\ .
\end{equation}

For each $j$, let $p_j$ satisfy $0 \le 1/p_j \le 1$, and define the vector $\vec p$
by
$$\vec p = (1/p_1,\dots,1/p_M)\ .$$
Define the constant $C(\vec p)$ by
$$C(\vec p) = \sup_{\vec f_1,\dots,\vec f_M \ne 0}\left\{
\frac{J(\vec f_1,\dots,\vec f_M) }{\prod_{j=1}^M|\vec f_j|_{p_j}}\right\}\ .$$

\medskip
\begin{thm}\label{thm6}
The function $\ln\left(C(\vec p)\right)$ is  convex on $[0,1]^M$.
\end{thm}

\medskip

We remark that one can extend the theorem to a version for multilinear functionals in $L^p$ 
spaces by the standard approximation argument with simple functions.

\medskip

\noindent{\bf Proof:}  Suppose that $\vec p$, $\vec q$ and $\vec r$ are vectors in $[0,1]^M$
such that for some $t$ with $0<t<1$,
$$\vec p = t\vec q +(1-t)\vec r\ .$$
Define the numbers $q_j$ and $r_j$ by $\vec q = (1/q_1,\dots,1/q_M)$ and $
\vec r = (1/r_1,\dots,1/r_M)$ so that $q_j$ is related to $\vec q$ and $r_j$ is related to $\vec r$ the same way $p_j$ is related to $\vec p$.  We must show that
$$C(\vec p) \le C(\vec q)^tC(\vec r)^{1-t}\ .$$

Towards this end, let $\vec f_1,\dots,\vec f_M$ be any $M$ non zero vectors. 
We can assume that 
$| \vec f_j| _{p_j}= 1$ for each $j$. Since the choice of the vectors is arbitrary apart from the normalization, it
suffices to show that
\begin{equation}\label{desire}
J(\vec f_1,\dots,\vec f_M) \le 
C(\vec q)^tC(\vec r)^{1-t}\ .
\end{equation}

Define
$$
\rho_{j,k} = |f_{j,k}|^{1/p}
$$
so that for each $j$ 
$$
\sum_{k=1}^N \rho_{j,k} = 1 \ .
$$
Also, define $\alpha_{j,k}$ by  $\alpha_{j,k} = f_{j,k}/|f_{j,k}|$ when $f_{j,k} \ne 0$, and 
$\alpha_{j,k} = 0$ otherwise.

Then, for each complex number $z$ define the vector
$ \vec \rho_j(z)$ by
\begin{equation}\label{60}
( \vec \rho_j(z))_k = \alpha_{j,k}\rho_{j,k}^{z/q_j+(1-z)/r_j} \ . 
\end{equation} 
Note that, for each $j$,
\begin{equation}\label{70}
 \vec \rho_j(t) = \vec f_j\ .
 \end{equation}
Moreover:

\smallskip
\noindent{\em (i)} The right side of (\ref{60}) is an entire function of $z$.

\smallskip
\noindent{\em (ii)}  Whenever $\Re (z) = 0$, then
$| \vec \rho_j(z) |_{r_j}=1$ 

\smallskip
\noindent{\em (iii)} 
Whenever $\Re(z)=1$ then $|  \vec \rho_j(z)|_{q_j}=1$.

Next we define
$$
G(z) = \left| \frac{ J(\vec \rho_1(z), \dots, \vec \rho_M(z))}{C(\vec q)^z C(\vec r)^{1-z}}\right|
$$
which is a subharmonic function of $z$. (By {\em (i)} and (\ref{sumform}), it is the absolute value of a sum of products of entire functions.)
 By {\em (ii)}, we know that on the line $\Re(z)=0, G(z) \le 1$,
and by  {\em (iii)}, we know that on the line $\Re(z)=1, G(z) \le 1$. Hence, by the maximum principle for
subharmonic functions, 
$G(t) \le 1$.  But by (\ref{70}), this yields (\ref{desire}). \qed

\end{document}